\documentclass[10pt,twoside,reqno]{amsart}
\usepackage{amssymb}
\usepackage{multirow}
\calclayout

\numberwithin{equation}{section}
\makeatletter

\renewcommand{\@secnumfont}{\bfseries}

\renewcommand{\section}{\@startsection{section}{1}%
  {0mm}{.7\linespacing\@plus\linespacing}{.5\linespacing}
  {\normalfont\bfseries\centering}}

\newcommand{\bibsection}{\@startsection{section}{1}%
  {0mm}{.7\linespacing\@plus\linespacing}{.5\linespacing}
  {\normalfont\scshape\centering}}

\renewcommand{\@biblabel}[1]{#1.}

\newtheorem{thm}{\bf Theorem}[section]
\newtheorem{lem}[thm]{\bf Lemma}
\newtheorem{cor}[thm]{\bf Corollary}

\begin{document}

\vspace{1.3cm}

\title { Some identities on derangement and degenerate derangement polynomials}

\author{Taekyun Kim}
\address{Department of Mathematics, College of Science Tianjin Polytechnic University, Tianjin 300160, China\newline
\indent Department of Mathematics, Kwangwoon University, Seoul 139-701, Republic
    of Korea}
\email{tkkim@kw.ac.kr}

\author{Dae San Kim}
\address{Department of Mathematics, Sogang University, Seoul 121-742, Republic
of Korea}
\email{dskim@sogang.ac.kr}

\begin{abstract}
In combinatorics, a derangement is a permutation that has no fixed points. The number of derangements of an $n$-element set is called the $n$-th derangement number.
In this paper, as natural companions to derangement numbers and degenerate versions of the companions we introduce derangement polynomials and degenerate derangement polynomials. We give some of their properties, recurrence relations and identities for those polynomials which are related to some special numbers and polynomials.
\end{abstract}
\maketitle
\bigskip
\medskip

\markboth{\centerline{\scriptsize  derangement and degenerate
derangement polynomials}} {\centerline{\scriptsize T. Kim, D. S.
Kim}}

\section{Introduction}

It is known that the Fubini polynomials are defined by the generating function
\begin{equation}\begin{split}\label{01}
\frac{1}{1-y(e^t-1)}=\sum_{n=0}^\infty F_n(y)\frac{t^n}{n!},\quad (\textnormal{see} \,\, [7,11]).
\end{split}\end{equation}
Thus, by \eqref{01}, we get
\begin{equation}\begin{split}\label{02}
F_n(y) = \sum_{k=0}^n S_2(n,k) k! y^k,\quad (\textnormal{see} \,\, [7,11]).
\end{split}\end{equation}
Here $S_2(n,k)$ is the Stirling number of the second kind which is defined by
\begin{equation}\begin{split}\label{03}
x^n = \sum_{l=0}^n S_2(n,l) (x)_l,\,\,(n \geq 0),
\end{split}\end{equation}
where $(x)_0=1, (x)_n = x(x-1) \cdots (x-n+1)$, $(n \geq 1)$.

As is well known, the Bell polynomials are given by the generating function as follows:
\begin{equation}\begin{split}\label{04}
e^{x(e^t-1)} = \sum_{n=0}^\infty           Bel_n(x)       \frac{t^n}{n!},\quad (\textnormal{see} \,\, [5,6,12]).
\end{split}\end{equation}
When $x=1$, $Bel_n = Bel_n(1)$ are called the Bell numbers. For $\lambda \in \mathbb{R}$, the partially degenerate Bell polynomials were introduced by Kim-Kim-Dolgy as
\begin{equation}\begin{split}\label{05}
e^{x\big( (1+\lambda t)^{\frac{1}{\lambda }}-1\big)} =\sum_{n=0}^\infty                  Bel_{n,\lambda }(x)\frac{t^n}{n!},\quad (\textnormal{see} \,\, [12]).
\end{split}\end{equation}

Note that $\lim_{\lambda \rightarrow 0} Bel_{n,\lambda }(x) = Bel_n(x)$, $(n \geq 0)$. When $x=1$, $Bel_{n,\lambda } = Bel_{l,\lambda }(1)$ are called the partially degenerate Bell numbers.

From \eqref{05}, we have
\begin{equation}\begin{split}\label{06}
Bel_{n,\lambda }(x) = \sum_{k=0}^n \sum_{m=0}^k S_2(k,m) S_1(n,k) \lambda ^{n-k} x^m,
\end{split}\end{equation}
where $S_1(n,k)$ is the Stirling number of the first kind given by
\begin{equation}\begin{split}\label{07}
(x)_n = \sum_{l=0}^n S_1(n,l) x^l,\,\,(n \geq 0),\quad (\textnormal{see} \,\, [8]).
\end{split}\end{equation}

In [1], L. Carlitz introduced the degenerate Bernoulli and Euler polynomials which are defined by
\begin{equation}\begin{split}\label{08}
\frac{t}{(1+\lambda t)^{\frac{1}{\lambda }}-1}(1+\lambda t)^{\frac{x}{\lambda }} = \sum_{n=0}^\infty           \beta_{n,\lambda }(x)       \frac{t^n}{n!},
\end{split}\end{equation}
and
\begin{equation}\begin{split}\label{09}
\frac{2}{(1+\lambda t)^{\frac{1}{\lambda }}+1}(1+\lambda t)^{\frac{x}{\lambda }} = \sum_{n=0}^\infty              \mathcal{E}_{n,\lambda }(x)    \frac{t^n}{n!}.
\end{split}\end{equation}
When $x=0$, $\beta_{n,\lambda }=\beta_{n,\lambda }(0)$, $\mathcal{E}_{n,\lambda }= \mathcal{E}_{n,\lambda }(0)$ are called the degenerate Bernoulli numbers and degenerate Euler numbers.

Recently, the degenerate Stirling numbers of the second kind are defined by
\begin{equation}\begin{split}\label{10}
S_{2,\lambda}(n+1,k) = kS_{2,\lambda }(n,k) + S_{2,\lambda }(n,k-1)-n\lambda S_{2,\lambda }(n,k),
\end{split}\end{equation}
where $n \geq 0$ (see [10]).

Note that $\lim_{\lambda \rightarrow 0}S_{2,\lambda }(n,k) = S_2(n,k)$. For $\lambda  \in \mathbb{R}$, the $\lambda$-analogue of falling factorial sequence is defined by
\begin{equation}\begin{split}\label{11}
(x)_{0,\lambda }=1,\,\,(x)_{n,\lambda }=x(x-\lambda )\cdots(x-(n-1)\lambda ),\,\,(n \geq 1),\quad (\textnormal{see} \,\, [6,8]).
\end{split}\end{equation}

Note that $\lim_{\lambda  \rightarrow 1} (x)_{n,\lambda } = (x)_n$, $(n\geq 0)$, (see [15]).

A derangement is a permutation with no fixed points. In other words, a derangement of a set leaves no elements in the original place. The number of derangements of a set of size $n$, denoted $d_n$, is called the $n$-th derangement number (see [9,13,14]).

For $n\geq 0$, it is well knwon that the recurrence relation of derangement numbers are given by
\begin{equation}\begin{split}\label{12}
d_n &= \sum_{k=0}^n {n \choose k} (n-k)! (-1)^k = n! \sum_{k=0}^n \frac{(-1)^k}{k!},\quad (\textnormal{see} \,\, [9]).
\end{split}\end{equation}

It is not difficult to show that
\begin{equation}\begin{split}\label{13}
\sum_{n=0}^\infty d_n \frac{t^n}{n!} = \frac{1}{1-t} e^{-t},\quad (\textnormal{see} \,\, [2,3,4,5,9]).
\end{split}\end{equation}

From \eqref{13}, we note that
\begin{equation}\begin{split}\label{14}
d_n = n \cdot d_{n-1}+(-1)^n,\,\,(n \geq 1),\quad (\textnormal{see} \,\, [9,13,14,16,17]).
\end{split}\end{equation}

and
\begin{equation}\begin{split}\label{15}
d_n = (n-1) (d_{n-1}+d_{n-2}),\,\,(n \geq 2).
\end{split}\end{equation}
In this paper, as natural companions to derangement numbers and degenerate versions of the companions we introduce derangement polynomials and degenerate derangement polynomials. We give some of their properties, recurrence relations and identities for those polynomials which are related to some special numbers and polynomials.

\section{Derangement polynomials}

Now, we define the derangement polynomials which are given by the generating function
\begin{equation}\begin{split}\label{16}
\frac{1}{1-xt}e^{-t} = \sum_{n=0}^\infty d_n(x) \frac{t^n}{n!}.
\end{split}\end{equation}

When $x=1$, $d_n(1)= d_n$ are the derangement numbers.

From \eqref{01}, we note that
\begin{equation}\begin{split}\label{17}
\frac{1}{1-yt} &= \sum_{m=0}^\infty F_m(y) \frac{1}{m!} \big( \log(1+t) \big)^m\\
&= \sum_{m=0}^\infty F_m(y) \sum_{n=m}^\infty S_1(n,m) \frac{t^n}{n!}\\
&= \sum_{n=0}^\infty \left( \sum_{m=0}^n F_m(y) S_1(n,m) \right)                 \frac{t^n}{n!}.
\end{split}\end{equation}

On the other hand,
\begin{equation}\begin{split}\label{18}
\frac{1}{1-yt} = \sum_{n=0}^\infty y^n n! \frac{t^n}{n!}.
\end{split}\end{equation}

Therefore, by \eqref{17} and \eqref{18}, we obtain the following lemma.

\begin{lem}
For $n \geq 0$, we have
\begin{equation*}\begin{split}
y^n = \frac{1}{n!}\sum_{m=0}^n F_m(y) S_1(n,m).
\end{split}\end{equation*}
\end{lem}

We observe that
\begin{equation}\begin{split}\label{19}
\frac{1}{1-yt}&= \left( \frac{1}{1-yt} e^{-t} \right) e^t=  \left( \sum_{l=0}^\infty d_l(y) \frac{t^l}{l!} \right) \left( \sum_{m=0}^\infty \frac{t^m}{m!} \right)\\
&= \sum_{n=0}^\infty   \left(          \sum_{l=0}^n {n \choose l} d_l(y)        \right)               \frac{t^n}{n!}.
\end{split}\end{equation}
From \eqref{17} and \eqref{19}, we obtain the following theorem.

\begin{thm}
For $n \geq 0$, we have
\begin{equation*}\begin{split}
 \sum_{l=0}^n {n \choose l} d_l(y)= \sum_{m=0}^n F_m(y) S_1(n,m).
\end{split}\end{equation*}
\end{thm}

By \eqref{16}, we get
\begin{equation}\begin{split}\label{20}
\sum_{n=0}^\infty d_n(x) \frac{t^n}{n!} &= \frac{1}{1-xt}e^{-t} = \left( \sum_{m=0}^\infty x^m t^m \right) \left( \sum_{k=0}^\infty \frac{(-1)^k}{k!} t^k \right)\\
&=\sum_{n=0}^\infty   \left(n! \sum_{k=0}^n \frac{(-1)^k}{k!} x^{n-k} \right) \frac{t^n}{n!}.
\end{split}\end{equation}
By comparing the coefficients on both sides of \eqref{20}, we obtain the following theorem.
\begin{thm}
For $n \geq 0$, we have
\begin{equation*}\begin{split}
d_n(x) =   n! \sum_{k=0}^n \frac{(-1)^k}{k!} x^{n-k}         .
\end{split}\end{equation*}
\end{thm}

From \eqref{16}, we have
\begin{equation}\begin{split}\label{21}
e^{-t} &= (1-xt) \sum_{n=0}^\infty d_n(x) \frac{t^n}{n!}\\
&=d_0(x) + \sum_{n=1}^\infty \left( d_n(x) - nx d_{n-1}(x) \right) \frac{t^n}{n!}.
\end{split}\end{equation}

On the other hand,
\begin{equation}\begin{split}\label{22}
e^{-t} = \sum_{n=0}^\infty (-1)^n \frac{t^n}{n!}.
\end{split}\end{equation}

Thus, by \eqref{21} and \eqref{22}, we get
\begin{equation}\begin{split}\label{23}
d_0(x) =1 ,\,\, d_n(x) = nx d_{n-1}(x) + (-1)^n,\,\,(n \geq 1).
\end{split}\end{equation}

From \eqref{23}, we note that
\begin{equation}\begin{split}\label{24}
d_n(x)&= (nx-1)d_{n-1}(x) + d_{n-1}(x) + (-1)^n\\
&= (nx-1)d_{n-1}(x) + (n-1)xd_{n-2}(x) + (-1)^{n-1} + (-1)^n\\
&= (nx-1) \left[ d_{n-1}(x) + d_{n-2}(x) \right]+(1-x)d_{n-2}(x),\,\,(n \geq 2).
\end{split}\end{equation}

Therefore, we obtain the following theorem.
\begin{thm}
For $n \geq 1$, we have
\begin{equation*}\begin{split}
d_n(x) = nxd_{n-1}(x) + (-1)^n.
\end{split}\end{equation*}
In particular, for $n \geq 2$, we have
\begin{equation*}\begin{split}
d_n(x) =(nx-1) \left[ d_{n-1}(x) + d_{n-2}(x) \right]+(1-x)d_{n-2}(x).
\end{split}\end{equation*}
\end{thm}

Replacing $t$ by $e^t-1$ in \eqref{16}, we get
\begin{equation}\begin{split}\label{25}
\frac{1}{1-x(e^t-1)}e^{-(e^t-1)} &= \sum_{m=0}^\infty d_m(x) \frac{1}{m!} (e^t-1)^m\\
 &= \sum_{m=0}^\infty d_m(x) \sum_{n=m}^\infty S_2(n,m) \frac{t^n}{n!}\\
 &= \sum_{n=0}^\infty  \left( \sum_{m=0}^n d_m(x) S_2(n,m)\right)\frac{t^n}{n!}.
\end{split}\end{equation}
By \eqref{25}, we see that
\begin{equation}\begin{split}\label{26}
\frac{1}{1-x(e^t-1)} &= e^{(e^t-1)} \sum_{k=0}^\infty \left(\sum_{m=0}^k d_m(x) S_2(k,m)\right) \frac{t^k}{k!}\\
&= \left( \sum_{l=0}^\infty Bel_l \frac{t^l}{l!} \right)\left(\sum_{k=0}^\infty  \left(\sum_{m=0}^k d_m(x) S_2(k,m)              \right) \frac{t^k}{k!}\right) \\
&= \sum_{n=0}^\infty  \left( \sum_{k=0}^n \sum_{m=0}^k {n \choose k} d_m(x) S_2(k,m) Bel_{n-k}                 \right)\frac{t^n}{n!}.
\end{split}\end{equation}

From \eqref{01}, we note that
\begin{equation}\begin{split}\label{27}
\frac{1}{1-x(e^t-1)} &= \sum_{n=0}^\infty F_n(x) \frac{t^n}{n!}.
\end{split}\end{equation}
Therefore, by \eqref{26} and \eqref{27}, we obtain the following theorem.

\begin{thm}
For $n \geq 0$, we have
\begin{equation*}\begin{split}
F_n(x) = \sum_{k=0}^n \sum_{m=0}^k {n \choose k} d_m(x) S_2(k,m) Bel_{n-k}           .
\end{split}\end{equation*}
\end{thm}

From \eqref{01}, we can derive the following equations \eqref{28}:
\begin{equation}\begin{split}\label{28}
\frac{1}{1-xt} e^{-t} &= \left( \sum_{k=0}^\infty \left( \sum_{m=0}^k F_m(x) S_1(k,m) \right) \frac{t^k}{k!} \right) e^{-t} \\
 &= \left( \sum_{k=0}^\infty \left( \sum_{m=0}^k F_m(x) S_1(k,m) \right) \frac{t^k}{k!} \right) \left( \sum_{l=0}^\infty \frac{(-1)^l}{l!} t^l \right)\\
 &= \sum_{n=0}^\infty   \left( \sum_{k=0}^n \sum_{m=0}^k {n \choose k} F_m(x) S_1(k,m) \frac{(-1)^{n-k}}{(n-k)!}                  \right)               \frac{t^n}{n!}.
\end{split}\end{equation}

On the other hand,
\begin{equation}\begin{split}\label{29}
\frac{1}{1-xt} e^{-t} &=\sum_{n=0}^\infty       d_n(x)           \frac{t^n}{n!}.
\end{split}\end{equation}

Therefore, by \eqref{28} and \eqref{29}, we obtain the following theorem.
\begin{thm}
For $n \geq 0$, we have
\begin{equation*}\begin{split}
d_n(x) = \sum_{k=0}^n \sum_{m=0}^k {n \choose k} F_m(x) S_1(k,m) \frac{(-1)^{n-k}}{(n-k)!}.
\end{split}\end{equation*}
\end{thm}

As is known, Bernoulli polynomials are defined by the generating function
\begin{equation}\begin{split}\label{30}
\frac{t}{e^t-1}e^{xt} = \sum_{n=0}^\infty B_n(x) \frac{t^n}{n!},\quad (\textnormal{see} \,\, [15]).
\end{split}\end{equation}
When $x=0$, $B_n=B_n(0)$ are Bernoulli numbers. By \eqref{30}, we easily get
\begin{equation}\begin{split}\label{31}
&\sum_{k=0}^{m-1} e^{kt} = \frac{1}{e^t-1} \big( e^{mt}-1 \big) = \frac{1}{t} \left\{ \frac{t}{e^t-1}e^{mt} - \frac{t}{e^t-1} \right\}\\
&= \sum_{n=0}^\infty \left( \frac{B_{n+1}(m)-B_{n+1}}{n+1} \right) \frac{t^n}{n!},\,\,(n \geq 1).
\end{split}\end{equation}

By Taylor expansion, we get
\begin{equation}\begin{split}\label{32}
\sum_{k=0}^{m-1} e^{kt} = \sum_{n=0}^\infty  \left( \sum_{k=0}^{m-1} k^n \right)                \frac{t^n}{n!}, \,\, (m \geq 1).
\end{split}\end{equation}

From \eqref{31} and \eqref{32}, we get
\begin{equation}\begin{split}\label{33}
\sum_{k=0}^{m-1} k^n = \frac{B_{n+1}(m)-B_{n+1}}{n+1} .
\end{split}\end{equation}

By Lemma 2.1, we easily get
\begin{equation}\begin{split}\label{34}
\sum_{k=0}^{m-1} k^n  =\frac{1}{n!}\sum_{k=0}^{m-1} \sum_{l=0}^n F_l(k) S_1(n,l).
\end{split}\end{equation}

Therefore, by Theorem 2.2, \eqref{33}, and \eqref{34}, we obtain the following theorem.

\begin{thm}
For $m \geq 1$ and $n \geq 0$, we have
\begin{equation*}\begin{split}
\frac{B_{n+1}(m)-B_{n+1}}{n+1}&=\frac{1}{n!}\sum_{k=0}^{m-1} \sum_{l=0}^n F_l(k) S_1(n,l)\\
&=\frac{1}{n!}\sum_{k=0}^{m-1} \sum_{l=0}^n {n \choose l} d_l(k).
\end{split}\end{equation*}
\end{thm}

\section{Degenerate derangement polynomials}

Here we consider the degenerate derangement polynomials which are given by
\begin{equation}\begin{split}\label{35}
\frac{1}{1-xt} (1-\lambda t)^{\frac{1}{\lambda }} = \sum_{n=0}^\infty                 d_{n,\lambda }(x) \frac{t^n}{n!},\,\,(\lambda  \in \mathbb{R}).
\end{split}\end{equation}

When $x=1$, $d_{n,\lambda }=d_{n,\lambda }(1)$ are called the degenerate derangement numbers.

From \eqref{35}, we note that
\begin{equation}\begin{split}\label{36}
(1-\lambda t)^{\frac{1}{\lambda }} &= \left(  \sum_{n=0}^\infty                 d_{n,\lambda }(x) \frac{t^n}{n!}               \right)(1-xt)\\
&= \sum_{n=0}^\infty       d_{n,\lambda }(x)           \frac{t^n}{n!}-\sum_{n=0}^\infty           x d_{n,\lambda }(x)       \frac{t^{n+1}}{n!}\\
&= d_{0,\lambda }(x) + \sum_{n=1}^\infty \left(      d_{n,\lambda }(x) - xnd_{n-1,\lambda }(x)            \right)                 \frac{t^n}{n!}.
\end{split}\end{equation}

On the other hand,
\begin{equation}\label{37}
(1-\lambda t)^{\frac{1}{\lambda }}= \sum_{m=0}^\infty { \frac{1}{\lambda } \choose m} (-\lambda )^m t^m = \sum_{m=0}^\infty (-1)^m (1)_{m,\lambda } \frac{t^m}{m!}.
\end{equation}

Therefore, by \eqref{36} and \eqref{37}, we obtain the following theorem.
\begin{thm}
For $n \geq 0$, we have
\begin{equation*}\begin{split}
d_{0,\lambda }(x) =1,\,\,d_{n,\lambda }(x) = nxd_{n-1,\lambda }(x) + (-1)^n (1)_{n,\lambda },\,\, (n \geq 1).
\end{split}\end{equation*}
\end{thm}

Note that $\lim_{\lambda \rightarrow 0} d_{n,\lambda }(x) = d_n(x)$, $\lim_{\lambda  \rightarrow 0}d_{n,\lambda }=d_n$, $(n \geq 0)$.

From \eqref{35}, we note that
\begin{equation}\begin{split}\label{38}
\sum_{n=0}^\infty         d_{n,\lambda }(x)         \frac{t^n}{n!} &=
\frac{1}{1-xt} (1-\lambda t)^{\frac{1}{\lambda }} = \left( \sum_{m=0}^\infty x^m t^m \right) \left( \sum_{k=0}^\infty (-1)^k (1)_{k,\lambda} \frac{t^k}{k!} \right)\\
&=\sum_{n=0}^\infty   \left( \sum_{k=0}^n \frac{(-1)^k}{k!} (1)_{k,\lambda } x^{n-k}                  \right)        t^n.
\end{split}\end{equation}

Comparing the coefficients on both sides of \eqref{38}, we obtain the following theorem.

\begin{thm}
For $n \geq 0$, we have
\begin{equation*}\begin{split}
d_{n,\lambda }(x) = n!\sum_{k=0}^n \frac{(-1)^k}{k!} (1)_{k,\lambda } x^{n-k}                .
\end{split}\end{equation*}
In particular, for $x=1$,
\begin{equation*}\begin{split}
d_{n,\lambda } = n!\sum_{k=0}^n \frac{(-1)^k}{k!} (1)_{k,\lambda } .
\end{split}\end{equation*}
\end{thm}

Now, we observe that
\begin{equation}\begin{split}\label{39}
\frac{1}{1-xt} &= \left( \frac{1}{1-xt} \right) (1-\lambda t)^{\frac{1}{\lambda }} \cdot (1-\lambda t)^{-\frac{1}{\lambda }} \\
&=\left( \sum_{l=0}^\infty d_{l,\lambda }(x) \frac{t^l}{l!} \right) \left( \sum_{m=0}^\infty { - \frac{1}{\lambda } \choose m} (-\lambda )^m t^m \right)\\
&=\left( \sum_{l=0}^\infty d_{l,\lambda }(x) \frac{t^l}{l!} \right) \left( \sum_{m=0}^\infty 1 (1+\lambda ) \cdots (1+(m-1)\lambda ) \frac{t^m}{m!} \right)\\
&= \sum_{n=0}^\infty  \left(\sum_{l=0}^n {n \choose l} d_{l,\lambda }(x) (1)_{n-l,-\lambda }           \right) \frac{t^n}{n!}.
\end{split}\end{equation}

On the other hand,
\begin{equation}\begin{split}\label{40}
\frac{1}{1-xt} = \sum_{n=0}^\infty x^n n! \frac{t^n}{n!}.
\end{split}\end{equation}

Therefore, by \eqref{39} and \eqref{40}, we obtain the following theorem.
\begin{thm}
For $n \geq 0$, we have
\begin{equation*}\begin{split}
x^n = \frac{1}{n!}  \sum_{l=0}^n {n \choose l} d_{l,\lambda }(x) (1)_{n-l,-\lambda }   .
\end{split}\end{equation*}
\end{thm}

From Theorem 3.1, we have
\begin{equation}\begin{split}\label{41}
d_{n,\lambda }(x) &= nxd_{n-1,\lambda }(x) + (-1)^n (1)_{n,\lambda } \\
&= (nx-1) d_{n-1,\lambda }(x) + d_{n-1,\lambda }(x) + (-1)^n (1)_{n,\lambda }\\
&= (nx-1) d_{n-1,\lambda }(x) + (n-1)x d_{n-2,\lambda }(x)\\
&\quad \quad+ (-1)^{n-1} (1)_{n-1,\lambda } + (-1)^n (1)_{n,\lambda }\\
&= (nx-1) \left[ d_{n-1,\lambda }(x) + d_{n-2,\lambda }(x) \right]\\
&\quad \quad+ (1-x)d_{n-2,\lambda }(x) + (-1)^{n-1}(1)_{n-1,\lambda }(n-1)\lambda,
\end{split}\end{equation}
where $n \geq 2$.

Therefore, by \eqref{41}, we obtain the following theorem.
\begin{thm}
For $n \geq 2$, we have
\begin{equation*}\begin{split}
&d_{n,\lambda }(x) = (nx-1) \left[ d_{n-1,\lambda }(x) + d_{n-2,\lambda }(x) \right]\\
&\quad \quad \quad \quad \quad+ (1-x)d_{n-2,\lambda }(x) + (-1)^{n-1}(1)_{n-1,\lambda }(n-1)\lambda.
\end{split}\end{equation*}
In particular, $x=1$,
\begin{equation*}\begin{split}
d_{n,\lambda } = (n-1) \left[ d_{n-1,\lambda } + d_{n-2,\lambda } \right] +\lambda (n-1) (-1)^{n-1}(1)_{n-1,\lambda }.
\end{split}\end{equation*}
\end{thm}
Note that
\begin{equation*}\begin{split}
d_n = \lim_{\lambda  \rightarrow 0} d_{n,\lambda } = (n-1) \left[ d_{n-1}+d_{n-2} \right]\,\,(n \geq 2).
\end{split}\end{equation*}

By using Taylor expansion, we get
\begin{equation}\begin{split}\label{42}
(1-\lambda t)^{\frac{1}{\lambda }}&= e^{\frac{1}{\lambda } \log(1-\lambda t)} = \sum_{m=0}^\infty \lambda ^{-m} \frac{1}{m!} \Big( \log(1-\lambda t)\Big)^m\\
&=\sum_{n=0}^\infty    \left(   \sum_{m=0}^n \lambda ^{n-m} (-1)^n S_1(n,m)               \right)              \frac{t^n}{n!}.
\end{split}\end{equation}

On the other hand,
\begin{equation}\begin{split}\label{43}
(1-\lambda t)^{\frac{1}{\lambda }} &= \frac{1}{1-xt} (1-\lambda t)^{\frac{1}{\lambda }} (1-xt)\\
&= \sum_{n=0}^\infty d_{n,\lambda }(x) \frac{t^n}{n!} - \sum_{n=1}^\infty                 nxd_{n-1,\lambda }(x) \frac{t^n}{n!}\\
&= d_{0,\lambda }(x) + \sum_{n=1}^\infty \left\{ d_{n,\lambda }(x) - nxd_{n-1,\lambda }(x) \right\}                 \frac{t^n}{n!}\\
&= 1+  \sum_{n=1}^\infty \left( d_{n,\lambda }(x) - nxd_{n-1,\lambda }(x) \right)                 \frac{t^n}{n!}\\
\end{split}\end{equation}

From \eqref{42} and \eqref{43}, we have
\begin{equation}\begin{split}\label{44}
(-1)^n \sum_{m=0}^n \lambda ^{n-m} S_1(n,m) =  d_{n,\lambda }(x) - nxd_{n-1,\lambda }(x) = (-1)^n (1)_{n,\lambda },\,\,(n \geq 1).
\end{split}\end{equation}

Therefore, by \eqref{44}, we obtain the following theorem.
\begin{thm}
For $n\geq 1$, we have
\begin{equation*}\begin{split}
\sum_{m=0}^n \lambda ^{n-m} S_1(n,m)              = (1)_{n,\lambda }.
\end{split}\end{equation*}
\end{thm}

By \eqref{13}, we get
\begin{equation}\begin{split}\label{45}
\frac{1}{(1+\lambda t)^{\frac{1}{\lambda }}+1} e^{(1+\lambda t)^{\frac{1}{\lambda }} } &= \sum_{m=0}^\infty (-1)^m d_m \frac{1}{m!} (1+\lambda t)^{\frac{m}{\lambda }}\\
&= \sum_{m=0}^\infty (-1)^m d_m \frac{1}{m!} \sum_{n=0}^\infty (m)_{n,\lambda } \frac{t^n}{n!}\\
&=\sum_{n=0}^\infty    \left( \sum_{m=0}^\infty (-1)^m d_m \frac{(m)_{n,\lambda }}{m!}                 \right)              \frac{t^n}{n!}.
\end{split}\end{equation}

On the other hand,
\begin{equation}\begin{split}\label{46}
\frac{1}{(1+\lambda t)^{\frac{1}{\lambda }}+1} e^{(1+\lambda t)^{\frac{1}{\lambda }} } &= \frac{e}{2} \frac{2}{(1+\lambda t)^{\frac{1}{\lambda }}+1} e^{(1+\lambda t)^{\frac{1}{\lambda }} -1} \\
&= \frac{e}{2} \left( \sum_{l=0}^\infty \mathcal{E}_{l,\lambda } \frac{t^l}{l!} \right) \left( \sum_{m=0}^\infty Bel_{m,\lambda } \frac{t^m}{m!} \right)\\
&= \frac{e}{2}\sum_{n=0}^\infty  \left(   \sum_{m=0}^n {n \choose m} Bel_{m,\lambda }\mathcal{E}_{n-m,\lambda }               \right)                \frac{t^n}{n!}.
\end{split}\end{equation}

Therefore, by \eqref{45} and \eqref{46}, we obtain the following theorem.
\begin{thm}
For $n \geq 0$, we have
\begin{equation*}\begin{split}
 \sum_{m=0}^n {n \choose m} Bel_{m,\lambda }\mathcal{E}_{n-m,\lambda }      =     \frac{2}{e} \sum_{m=0}^\infty (-1)^m d_m \frac{(m)_{n,\lambda }}{m!}         .
\end{split}\end{equation*}
\end{thm}

From \eqref{45}, we note that
\begin{equation}\begin{split}\label{47}
 e^{(1+\lambda t)^{\frac{1}{\lambda }}} &= \sum_{m=0}^\infty d_m \frac{(-1)^m}{m!} (1+\lambda t)^{\frac{m}{\lambda }} \Big(1+(1+\lambda t)^{\frac{1}{\lambda }}\Big)\\
 &= \sum_{m=0}^\infty d_m \frac{(-1)^m}{m!} (1+\lambda t)^{\frac{m}{\lambda }} + \sum_{m=0}^\infty d_m \frac{(-1)^m}{m!} (1+\lambda t)^{\frac{m+1}{\lambda }} \\
 &=\sum_{n=0}^\infty    \left\{ \sum_{m=0}^\infty d_m \frac{(-1)^m}{m!} \left( (m)_{n,\lambda } + (m+1)_{n,\lambda } \right)                   \right\}              \frac{t^n}{n!}.
\end{split}\end{equation}

On the other hand,
\begin{equation}\begin{split}\label{48}
 e^{(1+\lambda t)^{\frac{1}{\lambda }}} &= e \cdot  e^{(1+\lambda t)^{\frac{1}{\lambda }}-1} = e \sum_{k=0}^\infty \frac{1}{k!} \Big( (1+t)^{\frac{1}{\lambda }}-1 \Big)^k \\
 &=e \sum_{k=0}^\infty \sum_{n=k}^\infty S_{2,\lambda }(n,k) \frac{t^n}{n!} = e \sum_{n=0}^\infty    \left( \sum_{k=0}^n S_{2,\lambda }(n,k)                 \right)              \frac{t^n}{n!}.
\end{split}\end{equation}

Therefore, by \eqref{47} and \eqref{48}, we obtain the following theorem.
\begin{thm}
For $n \geq 0$, we have
\begin{equation*}\begin{split}
\sum_{m=0}^n S_{2,\lambda }(n,m) = \frac{1}{e} \sum_{m=0}^\infty d_m \frac{(-1)^m}{m!} \left( (m)_{n,\lambda } + (m+1)_{n,\lambda } \right).
\end{split}\end{equation*}
\end{thm}

Indeed,
\begin{equation}\begin{split}\label{49}
\sum_{n=0}^\infty    Bel_{n,\lambda }              \frac{t^n}{n!} &=
e^{\Big((1+\lambda t)^{\frac{1}{\lambda }}-1\Big)} = \sum_{m=0}^\infty \frac{1}{m!} \Big( (1+\lambda t)^{\frac{1}{\lambda }}-1 \Big)^m \\
&= \sum_{m=0}^\infty \sum_{n=m}^\infty S_{2,\lambda }(n,m) \frac{t^n}{n!}= \sum_{n=0}^\infty      \left(     \sum_{m=0}^n S_{2,\lambda }(n,m)             \right)            \frac{t^n}{n!}.
\end{split}\end{equation}

Thus, by \eqref{49}, we get
\begin{equation}\begin{split}\label{50}
Bel_{n,\lambda } = \sum_{m=0}^n S_{2,\lambda }(n,m),\,\,(n \geq 0).
\end{split}\end{equation}

Therefore, by \eqref{50}, we obtain the following corollary.
\begin{cor}For $n \geq 0$, we have
\begin{equation*}\begin{split}
Bel_{n,\lambda } = \frac{1}{e} \sum_{m=0}^\infty d_m \frac{(-1)^m}{m!} \left( (m)_{n,\lambda } + (m+1)_{n,\lambda } \right).
\end{split}\end{equation*}
\end{cor}

\end{document}